# Price-Based Market Clearing with V2G Integration Using Generalized Benders Decomposition

Reza Jamalzadeh, *Member, IEEE*, Sajjad Abedi, *Member, IEEE,* Masoud Rashidinejad, *Senior Member*, *IEEE,* and Mingguo Hong, *Member, IEEE*

*Abstract*— Currently, most ISOs adopt offer cost minimization (OCM) auction mechanism which minimizes the total offer cost, and then, a settlement rule based on either locational marginal prices (LMPs) or market clearing price (MCP) is used to determine the payments to the committed units, which is not compatible with the auction mechanism because the minimized cost is different from the payment cost calculated by the settlement rule. This inconsistency can drastically increase the payment cost. On the other hand, payment cost minimization (PCM) auction mechanism eliminates this inconsistency; however, PCM problem is a nonlinear self-referring NP-hard problem which poses grand computational burden. In this paper, a mixed-integer nonlinear programing (MINLP) formulation of PCM problem are presented to address additional complexity of fast-growing penetration of Vehicle-to-Grid (V2G) in the price-based market clearing problem, and a solution method based on the generalized benders decomposition (GBD) is then proposed to solve the V2G-integrated PCM problem, and its favorable performance in terms of convergence and computational efficiency is demonstrated using case studies. The proposed GBD-based method can handle scaled-up models with the increased number of decision variables and constraints which facilitates the use of PCM mechanism in the market clearing of large-scale power systems. The impact of using V2G technologies on the OCM and PCM mechanisms in terms of MCPs and payments is also investigated, and by using numerical results, the performances of these two mechanisms are compared.

## Nomenclature

*Parameters*

| | |
|---|---|
| $B_i(t)$ | Single block bid price of unit *i* at time period *t*. |
| $C_i^{NL}$ | No load cost of generating unit *i*. |
| $D(t)$ | Total system load demand at time period *t*. |
| $E_v^{min}, E_v^{max}$ | Minimum/maximum energy stored in batteries of PEV fleet *v*. |
| $NT_i$ | Number of intervals of the stair-wise startup cost function of unit *i*. |
| $P_i^{max}(t), P_i^{min}(t)$ | Maximum and minimum power offered by unit *i* at time period *t*. |
| $P_v^{CH,max}, P_v^{CH,min}$ | Maximum/minimum charging capacity of PEV fleet *v*. |
| $P_v^{DSCH,max}, P_v^{DSCH,min}$ | Maximum/minimum discharging capacity of PEV fleet *v*. |
| $RU_i, RD_i$ | Ramp up and ramp down rates for unit *i*. |
| $SU_i^{t_i^{off}}$ | Cost of the interval $t_i^{off}$ (offline time) of the stair-wise startup cost function of unit *i*. |
| $SD_i$ | Shutdown cost of unit *i*. |
| $t_i^{off}$ | Number of periods unit *i* has been offline prior to the startup |
| $T$ | Number of periods of the time span. |
| $\eta_v$ | Charging/discharging cycle efficiency of the PEV fleet *v*. |
| $\lambda, \mu$ | Lagrangian multipliers. |

*Variables*

| | |
|---|---|
| $E_v(t)$ | Available energy in batteries of fleet *v* at time *t*. |
| $MCP(t)$ | MCP at time *t*. |
| $p_i(t)$ | Power output of unit *i* at time *t*. |
| $p_v(t)$ | Power output of PEV fleet *v* at time *t*. |
| $sc_i^u(t)$ | Startup cost of unit *i* at time *t*. |
| $sc_i^d(t)$ | Shutdown cost of unit *i* at time *t*. |
| $u_i(t)$ | Commitment status of unit *i* at time *t*. |
| $u_v^{CH}(t), u_v^{DSCH}(t)$ | Charging/discharging status of PEV fleet *v* at time *t*. |

*Set Indices*

| | |
|---|---|
| *i* | Index for generating units. |
| *t* | Index for study time interval. |
| *v* | Index for PEV fleets. |

## I. Introduction

### A. Background and Motivation

In electricity markets (e.g., the day-ahead markets), based on offers received from market participants (i.e., energy offers from producers and energy bids from consumers), independent system operators (ISOs) use a clearing algorithm to determine the market-clearing price, the power productions, and the consumption level of consumers in every period of time. Generally, there are two main clearing mechanisms: First, offer cost minimization (OCM) mechanism which is used to select offers in a way that the total bid cost is minimized, and second, payment cost minimization (PCM) which is used to select offers for minimizing total actual payments to the accepted bidders [1-7]. After using clearing algorithm by ISO, in order to determine the payments to the selected bidders, a settlement rule (e.g., pay-as-bid pricing, uniform pricing) should be used [1-7]. In pay-as-bid pricing, since each accepted bidder is paid at its

Reza Jamalzadeh and Mingguo Hong are with Department of Electrical Engineering and Computer Science, Case Western Reserve University, OH, USA. (e-mail: rxj171@case.edu and mxh543@case.edu). Sajjad Abedi is with School of Mechanical Engineering, Purdue University, IN, USA. (email: sabedi@purdue.edu). Masoud Rashidinejad is with Department of Electrical Engineering, Bahonar University, Iran. (e-mail: mrashidi@uk.ac.ir).

offer, the payment cost is the same as the offer cost so that there is no difference between two auction mechanisms. However, when pay-as-market clearing price (MCP) or pay-as-locational marginal price (LMP) is utilized as the settlement rule, the payment cost would be different from the offer cost; therefore, PCM and OCM auctions may provide different clearing solutions which might result in different market equilibrium.

Currently, the majority of market operators (e.g., NY-ISO, ISO-NE, ERCOT, PJM, MISO) adopt OCM auction mechanism in their markets. One of the justifications is that the OCM model is simpler to solve as compared to the PCM, and straightforward solution approaches have been developed to minimize the total production cost (including units' offer costs and fixed costs) and maximize social welfare. Afterwards, pay-as-MCP mechanism or pay-as-LMP mechanism is used as a settlement rule which is not compatible with the auction mechanism utilized because the minimized cost is different from the payment cost. This inconsistency between the clearing algorithm (i.e., OCM) and the settlement rule (i.e., pay-as-MCP or pay-as-LMP) can drastically increase the payment cost. By using PCM auction mechanism, this inconsistency is effectively eliminated, and considerable reduction in payments in comparison with payments of OCM mechanism is achieved [1-10]. Literature has shown that for the same set of bids, PCM leads to reduce consumer payments. The market participants, however, may bid differently under the two auction mechanisms. In [11], the supplier's strategic behaviors are investigated in a simplified day-ahead energy market under the two auctions, and it is concluded that PCM still leads to significant reductions in payments even with strategic bidding. In [12], it is shown that the sensitivities of LMPs with respect to the system uncertainties under the PCM mechanism are lower comparing to those under the OCM mechanism. This demonstrates yet another significant advantage of the PCM over the OCM mechanism.

The PCM problem has recently received considerable attention due to the open challenges in both modeling and solution algorithm. The PCM problem is an NP-hard problem [3, 7]. Also because the market prices (either MCP or LMP) are present in the objective function, the PCM problem is a self-referring optimization problem [12] and suffers from added complexity and computational burden. Therefore, inefficient solution algorithms may lead to prolonged computational time or even failure of convergence. To address these challenges, this paper presents an efficient solution algorithm based on the Generalized Benders Decomposition (GBD).

*B. Related Work on PCM*

The existing OCM auction in electricity markets is similar to the unit commitment problem in centralized market operations and solution methods for solving OCM auction abound. The solution methods for solving PCM problem, however, are limited and mostly inefficient [3, 4]. In [8], a solution method based on forward dynamic programming was presented to solve a simple PCM problem, but the author acknowledges that the method is not suited for large-scale problems ( due to "curse of dimensionality"). Reference [5] proposed a graph search algorithm to solve a simple PCM problem by assuming simple bids with price-quantity curve. But because of the complexity of the method, it is not suitable for solving large-scale problems either. In addition, the MCPs in this study are loosely defined as the maximum of amortized bid costs. In [13], the genetic algorithm was used to solve the PCM problem that also suffers from low computational efficiency. References [3, 7, 9] presented solution methods based on the augmented Lagrangian relaxation and surrogate optimization. But the proposed methods cannot guarantee solution feasibility [4]. In [10], a solution method based on bilevel programing was proposed and in [4], we presented a mixed-integer linear programming (MILP) formulation of PCM problem. The solution methods presented in [4] and [10] can merely be used to solve PCM problems with inelastic demands, but they cannot deal with the nonlinearities posed by energy storage systems, demand response (or elastic demand), and V2G. In this paper we extend the formulation presented in [4] to address this challenge, and solve the resultant mixed integer nonlinear programing (MINLP) formulation of PCM through GBD.

It is observed that some studies in the literature [1-5, 8, 9] solved the PCM problem without considering the transmission constraints, while some others did [7, 10, 12]. Without loss of generality, transmission constraints are not considered in this study and the uniform price settlement mechanism based on MCP is used. As necessary, transmission constraints can be added to the model and the LMPs can be used as market clearing prices.

*C. Integration of V2G in PCM*

Plug-in electric vehicles (PEVs) as the portable source of electricity storages have undeniable benefits through intelligent charging and discharging scheme in a smart grid environment [14-20]. The techno-economic advantages of PEVs include flattening load curve and minimizing load curtailment by discharging PEVs in peak time and charging at off-peak periods, and providing frequency regulation and spinning reserves in fast response, etc. Due to the significant benefits, there is a pressing need to study market operations with a high level of V2G penetration. This paper has studied the impact of V2G on OCM and PCM in terms of both MCPs and payments for the first time.

*D. Main Contributions of This Work*

In our proposed GBD-based solution algorithm, the MINLP model formulation of PCM splits into 1) feasibility and optimality subproblems that are mixed integer linear programing (MILP) problems, and 2) a master problem which is a linear programing (LP) problem. The case studies and numerical comparison with other solution methods confirm the solution performance of the proposed solution algorithm, in terms of both computational convergence and speed. The proposed GBD-based method can handle scaled-up models with considerations for V2G, demand response (or elastic demand), and energy storage systems with the increased number of decision variables and constraints which facilitates the use of PCM mechanism in the market clearing of large-scale power systems. Meanwhile, this work has modeled V2G in the PCM auction mechanisms for the first time. After analyzing the impacts on MCPs and payments, our study demonstrates further benefits of the PCM over the OCM mechanism with V2G integration.

The remainder of the paper is organized as follows. Section II presents the MINLP formulation of PCM-V2G problem; for comparison, the MILP formulation of the OCM-V2G problem is also described. Section III presents the GBD-based solution method for the PCM-V2G problem. Case studies and numerical results are illustrated in section IV where the benefits of V2G on both the PCM and OCM mechanisms are compared while the GBD-based algorithm solution performance is demonstrated. Finally, conclusions are presented in section V.

## II. OCM AND PCM FORMULATIONS

In this section, the formulations of OCM and PCM auction mechanisms incorporating V2G are presented. It is assumed that system load demand is fixed and transmission constraints are not considered without loss of generality.

### A. The OCM-V2G Problem

Generally, the objective function of the OCM-V2G problem is formulated as the minimization of the total production cost as follows:

$$\min \sum_t \sum_i \left( O_i(p_i(t),t) \cdot p_i(t) + C_i(t) \right) \quad (1)$$

where $O_i(p_i(t),t)$ is the offer price function of unit $i$ in terms of $p_i(t)$ at time period $t$, and $C_i(t)$ represents the fixed costs associated with unit $i$ at time $t$. By considering the typical single block offer curve, $O_i(p_i(t),t)$ in (1) as a generic function is replaced with $B_i(t)$, and the objective function of OCM problem becomes:

$$\min \left\{ \sum_{t=1}^{T} \sum_{i=1}^{I} \left( B_i(t) \cdot p_i(t) + sc_i^u(t) + sc_i^d(t) + C_i^{NL} \cdot u_i(t) \right) \right\} \quad (2)$$

where $B_i(t)$ is a constant parameter for unit $i$ at time $t$, and $sc_i^u(t)$ and $sc_i^d(t)$ are each the startup cost and shutdown cost variables associated with unit $i$, and $C_i^{NL}$ is the no load cost of unit $i$. The OCM-V2G problem is subject to the following constraints.

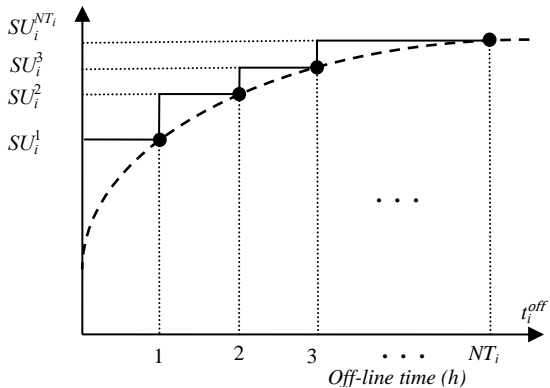

Fig. 1. Exponential and stair-wise startup cost functions.

*a.* Startup cost constraints:
A typical exponential startup cost function in terms of offline time is shown in Fig. 1 by the dashed line. To model the problem as a MILP formulation, the stair-wise startup cost function (solid line in Fig.1) is modeled by using the following two constraints.

$$sc_i^u(t) \geq SU_i^{t_i^{off}}(t) \left( u_i(t) - \sum_{n=1}^{t_i^{off}} u_i(t-n) \right), \forall i, \forall t, \forall t_i^{off} = 1,...,NT_i \quad (3)$$

$$sc_i^u(t) \geq 0 \;, \forall i, \forall t \quad (4)$$

where $SU_i^{t_i^{off}}$ is the startup cost of the interval $t_i^{off}$ of the stair-wise startup cost function of unit $i$ (Fig. 1).

*b.* Shut down cost constraints:

$$sc_i^d(t) \geq SD_i \cdot \left( u_i(t-1) - u_i(t) \right) \;, \forall i, \forall t \quad (5)$$

$$sc_i^d(t) \geq 0 \;, \forall i, \forall t \quad (6)$$

where constant parameter $SD_i$ is the shutdown cost of unit $i$.

*c.* Power balance equation constraints:

$$D(t) = \sum_i p_i(t) + \sum_v p_v(t) \;, \forall t \quad (7)$$

where $p_v(t)$ is the output power of PEV fleet $v$ at time $t$. Variable $p_v(t)$ is positive when the aggregated PEV fleet is discharging, and negative when the aggregated PEV fleet is charging.

*d.* Constraints for the power output of generating units:
$$P_i^{\min}(t) \cdot u_i(t) \leq p_i(t) \leq P_i^{\max}(t) \cdot u_i(t), \forall i, \forall t \quad (8)$$

*e.* Active power ramp constraints:
$$p_i(t) - p_i(t-1) \leq u_i(t-1) \cdot RU_i + [1 - u_i(t-1)] \cdot P_i^{\max}, \forall i, \forall t \quad (9)$$
$$p_i(t-1) - p_i(t) \leq u_i(t) \cdot RD_i + [1 - u_i(t)] \cdot P_i^{\max}, \forall i, \forall t \quad (10)$$

*f.* Constraints for charging and discharging power of the PEV fleets:
$$p_v(t) \leq u_v^{DSCH}(t) \cdot P_v^{DSCH,\max} - u_v^{CH}(t) \cdot P_v^{CH,\min}, \forall t, \forall v \quad (11)$$
$$p_v(t) \geq u_v^{DSCH}(t) \cdot P_v^{DSCH,\min} - u_v^{CH}(t) \cdot P_v^{CH,\max}, \forall t, \forall v \quad (12)$$
$$u_v^{CH}(t) + u_v^{DSCH}(t) \leq 1, \forall t, \forall v \quad (13)$$

where $u_v^{CH}(t)$ and $u_v^{DSCH}(t)$ are binary variables representing the charging and discharging status of PEV fleet $v$ at time $t$, respectively.

*g.* The energy equation constraint of each PEV fleet:
$$E_v(t) = E_v(t-1) - \eta_v \cdot p_v(t), \forall t, \forall v \quad (14)$$

where $E_v(t)$ is the available energy in batteries of fleet $v$ at time $t$ and is calculated based the available energy at time ($t$-1) and the hourly power output of PEV fleet $v$ factored by the charging cycle efficiency of the PEV fleet $v$ ($\eta_v$).

*h.* Constraints for energy capacity limits of the PEV fleets:
$$E_v^{\min} \leq E_v(t) \leq E_v^{\max}, \forall t, \forall v \quad (15)$$
$$E_v(T) = \sigma_v(T), \forall v \quad (16)$$

where $E_v^{\min}(t)$ and $E_v^{\max}(t)$ are minimum and maximum energy stored in batteries of PEV fleet $v$, and $\sigma_v(T)$ is the targeted energy level of PEV fleet $v$ at the end of the study period.

Other generator constraints can be also incorporated, such as minimum run time [21]. There are numerous solution methods for the OCM-V2G problem. In this study, the OCM problem as defined by (2)-(16) is solved as a MILP problem using a readily available commercial MILP solver such as CPLEX [22].

After clearing the market, most ISOs use pay-as-MCP or pay-as-LMP as the settlement rule. When transmission constraints are not considered, the uniform price settlement mechanism based on MCP is used. As a result, the MCP is defined as the highest offer accepted during each time period $t$.

$$MCP(t) = \max \{O_i(P_i(t),t) \cdot u_i(t)\} \quad (17)$$

The MCP as determined by (17) is merely used for energy payment; the committed generators can be paid separately to cover startup costs and shut down costs.

### B. The PCM-V2G Problem

In the PCM-V2G problem without transmission constraints, a uniform market clearing price, i.e., MCP is being determined to clear the market. The objective function of the PCM problem directly minimizes the procurement cost and can be generally formulated as the following:

$$\min \sum_t \sum_i \left(MCP(t) \cdot p_i(t) + C_i(t)\right) \quad (18)$$

where $C_i(t)$ represents the fixed costs associated with unit $i$ at time $t$. In (18), $MCP(t)$ is equal to the highest offer accepted by ISO. In the PCM-V2G problem formulation, $MCP(t)$ is defined as a variable. With specific fixed cost, the objective function of PCM-V2G is presented as:

$$\min \sum_t \sum_i \left(MCP(t) \cdot p_i(t) + sc_i^u(t) + sc_i^d(t) + C_i^{NL} \cdot u_i(t)\right) \quad (19)$$

The following constraint ensures that the value of $MCP(t)$ is the maximum accepted offer among all generation offers $B_i(t)$ [4].

$$MCP(t) \geq B_i(t) \cdot u_i(t), \forall i, \forall t \quad (20)$$

Under constraint (20), continuous variables $MCP(t)$ take on discrete values. The PCM-V2G problem is defined by objective function (19), and constraints (3)-(16), and (20).

## III. PROPOSED SOLUTION METHOD USING GENERALIZED BENDERS DECOMPOSITION

The nonlinearity of objective functions presented in (19) lines in the product of two variables in the first part of the objective function. This problem is a NP-hard problem; if the problem dimension becomes larger with more variables, solving the problem gets more difficult or even impossible, such as for a large-scale system with various DG sources. In this study, GBD [23] is utilized to split the problem into subproblems which are MILP problems, and master problem which is a LP problem.

### A. Formulation of the optimality subproblem

To apply the GBD method to the problem, the decision variables are split into two groups $X$ and $Y$.

$$Y = \{MCP(t)\}$$
$$X = \{p_i(t), u_i(t), sc_i^u(t), sc_i^d(t), p_v(t), u_v^{CH}(t), u_v^{DSCH}(t), E_v(t)\} \quad (21)$$

where $Y$ is the set of complicating variables. Then, the objective function and constraints of PCM-V2G can be written as:

$$\min_{x \in X, y \in Y} \sum_t \sum_i \left(MCP(t) \cdot p_i(t) + sc_i^u(t) + sc_i^d(t) + C_i^{NL} \cdot u_i(t)\right) \quad (22)$$

$$MCP(t) \geq B_i(t) \cdot u_i(t) \quad (22\text{-}1)$$
$$g(X) \leq 0 \quad (22\text{-}2)$$

where $g(X) \leq 0$ corresponds to the constraints (3)-(16). By fixing $y_i = \hat{y}_i \in Y$ and for $x_i \in X$, the optimality subproblem is formulated as:

$$\min_{x \in X} \sum_t \sum_i \left(\hat{MCP}(t) \cdot p_i(t) + sc_i^u(t) + sc_i^d(t) + C_i^{NL} \cdot u_i(t)\right) \quad (23)$$

$$\hat{MCP}(t) \geq B_i(t) \cdot u_i(t) \qquad \lambda_i(t) \quad (23\text{-}1)$$
$$g(X) \leq 0 \qquad \mu \quad (23\text{-}2)$$

where $\lambda_i(t)$ is the lagrangian multiplayer corresponding to constraint (23-1) for unit $i$ at time period $t$, and $\mu$ is the vector of lagrangian multiplayer corresponding to constraint (23-2). Now, if optimality subproblem (23) is feasible, its optimal solution is denoted by $\hat{X}$ with optimal multiplier vectors $\hat{\mu}$ and $\hat{\lambda}$. The solution of (23) provides an upper bound (UBD) for the solution of the main problem (22).

$$UBD = \hat{Z}_{Sub} \quad (24)$$

where UBD is the upper bound of the main problem (22) and $\hat{Z}_{sub}$ is the optimal objective function value of the optimality subproblem (23).

### B. Formulation of the feasibility subproblem

If the optimality subproblem (23) is infeasible, the feasibility subproblem is formulated as:

$$\min_{x \in X, \alpha_i(t)} w = \sum_t \sum_i \alpha_i(t) \quad (25)$$

$$\hat{MCP}(t) + \alpha_i(t) \geq B_i(t) \cdot u_i(t) \qquad \lambda_i(t) \quad (25\text{-}1)$$
$$g(X) \leq 0 \qquad \mu \quad (25\text{-}2)$$

where $\alpha_i(t)$ is the slack variable corresponding to unit $i$ at time $t$, and $g(X) \leq 0$ corresponds to constraints (3)-(16).

### C. Formulation of the master problem

Each time the optimality or feasibility subproblem is solved, an optimality cut or a feasibility cut is generated and enforced for the following master problem:

$$\min_{y \in Y, \eta} \eta \quad (26)$$

s.t.

Optimality cut: $\eta \geq L^*(Y, \hat{\mu}, \hat{\lambda}) \quad (26\text{-}1)$

Feasibility cut: $L_*(Y, \hat{\mu}, \hat{\lambda}) \leq 0 \quad (26\text{-}2)$

with:

$$L^*(Y, \mu, \lambda) = \inf_X \left\{ \begin{array}{l} \sum_t \left(MCP(t) \cdot \sum_i p_i(t)\right) + \sum_t \sum_i \left(\begin{array}{l} sc_i^u(t) + sc_i^d(t) + \\ C_i^{NL} \cdot u_i(t) \end{array}\right) + \\ \sum_t \sum_i \left(\lambda_i(t) \cdot \left(B_i(t) \cdot u_i(t) - MCP(t)\right)\right) + \mu^T \cdot g(X) \end{array} \right\}$$

$$(26\text{-}3)$$

$$L_*(Y,\mu,\lambda) = 1^T\alpha + \inf_X \left\{ \begin{array}{l} \sum_t\sum_i \lambda_i(t)\cdot(B_i(t)\cdot u_i(t) - MCP(t) - \alpha_i(t)) + \\ \mu^T \cdot g(X) \end{array} \right\}$$
(26-4)

where multipliers $\hat{\mu}$ and $\hat{\lambda}$ in (26-1) and (26-2) are iteratively provided from the solutions of optimality and feasibility subproblems. The explicit expressions of function $L^*$ presented in (26-3) can't be easily obtained since variables $x$ and $y$ are not separable inside the infimum term which satisfy neither the $P$ nor the $P'$ properties [24]. However, as often practiced in other engineering fields [24], the primal optimal solution $\hat{x}$ of the optimality subproblem is utilized in this study to approximate $L^*$, and thus, the optimality cut (26-1) can be simplified as follows:

$$\eta \geq UBD + \sum_t\sum_i \pi_{i,t}\left(MCP(t) - M\hat{C}P(t)\right) \quad (27)$$

where $\pi_{i,t}=\hat{\lambda}_{i,t}+\hat{p}_i(t)$ in which $\pi_{i,t}$, $\hat{\lambda}_{i,t}$, and $\hat{p}_i(t)$ are each the optimality multiplayer, optimal lagrangian multiplayer corresponding to (23-1), and optimal generation schedule of unit $i$ at time $t$ obtained from the optimality subproblem (23) for fixed $y=\hat{y}$. Likewise, feasibility cuts can be formulated for each time period $1\leq t \leq T$ as follows:

$$\hat{w}_t + \sum_i \hat{\lambda}_{i,t}\cdot\left(MCP(t) - M\hat{C}P(t)\right) \leq 0, \forall t \quad (28)$$

where $\hat{w}_t=\sum_i\hat{\alpha}_i(t)$ in which $\hat{\alpha}_i(t)$ is the slack variable corresponding to unit $i$ at time period $t$ achieved from the feasibility subproblem (25) for fixed $y=\hat{y}=M\hat{C}P(t)$. In (28), $\hat{\lambda}_{i,t}$ is the optimal lagrangian multiplayer corresponding to (25-1) obtained from the feasibility subproblem (25) for fixed $y=\hat{y}=M\hat{C}P(t)$. The solution of the master problem (26) provides a lower bound (*LBD*) for the solution of main problem (22).

$$LBD = \hat{\eta} \quad (29)$$

where $\hat{\eta}$ is the optimal objective function value of the master problem (26).

*D. Summary of the GBD-based solution algorithm*

The GBD-based solution algorithm is started using an initial guess $y=\hat{y}=M\hat{C}P(t)$. For example, the initial point can be the highest bid price of generation units, i.e., $y_0(t)=M\hat{C}P(t)=\max_i\{B_i(t)\}$. As shown in [4], the solution of OCM problem can also be used as an initial value for the PCM problem. Then, optimality subproblem (23) is solved for the initial guess $\hat{y}=y_0$, and *UBD* for the main problem (22) is obtained through (24), and an optimality cut for the master problem is constructed. If optimality subproblem (23) is infeasible, the feasibility subproblem (25) is solved and a feasibility cut for the master problem is constructed. After solving the subproblem, the master problem is solved with the new optimality or feasibility cut, and *LBD* for the main problem (22) is obtained through (29), and the updated solution $\hat{y}=y_1$ for subproblems are achieved. The solution procedure continues, and optimality subproblem or feasibility subproblem is constructed with the updated solution $\hat{y}=y_1$. The above iterative process continues until the *UBD* and *LBD* converge within a predefined tolerance threshold.

## IV. NUMERICAL CASE STUDIES

In this study, a 24-hour problem with 10 generating units is considered. The system data, market offers, and hourly system load demands have been extracted from [21, 25]. In this case study, startup cost functions are modeled by a two interval stair wise linear function.

$$SU_i^{t_{off}} = \begin{cases} HSC_i & \forall t_{off} \leq T_i^{int} \\ CSC_i & \forall t_{off} \geq T_i^{int} \end{cases} \quad (27)$$

where $HSC_i$ is the hot startup cost, and $CSC_i$ is the cold startup cost of unit $i$; therefore, startup cost of first and second block of stair wise curve in Fig. 1 are equal to $HSC_i$ and $CSC_i$, respectively. In fact, if offline time of unit $i$ is lower than $T_i^{int}$, the startup cost is equal to $HSC_i$, and startup cost is equal to $CSC_i$, otherwise.

Also in this study, the 100000 PEVs are assumed for simulation. The probabilistic model of the attendance of PEVs in the PEV fleets for the hours of a day is extracted from [26], where the data is obtained from a parking lot from the city of Livermore, CA [27]. The maximum number of PEVs with charging/discharging status at each hour is presumed to be 10 percent of total PEVs in PEV fleets. In addition, the following parameters are considered for each vehicle in the PEV fleets: battery capacity of each PEV is 15kWh, state of charge of each PEV is 50%, and charging/discharging cycle efficiency of each PEV is 85%.

In model implementation, GAMS [28] with the CPLEX solver [22] is used to implement and solve both the OCM-V2G and the GBD-based PCM-V2G problems. The algorithm solutions are performed on desktop computers with Intel Core i3 CPU 530 @3 GHz and 3 GB RAM.

For comparative analysis, the OCM and PCM problems are first solved as MIP problems without considering V2G. The load demand profile, MCP as well as the payment cost of both PCM and OCM as results of the solution are presented in Fig. 2. As shown in Fig. 2, the PCM problem solves to lower MCPs and payment costs in the valley of the load profile, as compared to the OCM problem. Also as shown in Table I, the average MCP (62.792 $/MWh) and the total payment cost ($1816480) of the PCM problem are lower than those of the OCM (63.125 $/MWh and $1823570, respectively). Therefore without the V2G, the PCM mechanism achieves a total saving of $7090 (0.38%) over 24-hour study period as compared to the OCM mechanism.

The case study is then examined with considering V2G in the system. To solve the PCM-V2G problem, the proposed GBD-based method is utilized to split the main problem, which is a self-referring and NP-hard problem, into MILP optimality and feasibility subproblems as well as LP master problem. Solver CPLEX is used to solve the resulted subproblems and master problem in each iteration. The solution algorithm starts with an initial guess for $MCP(t)$. Based on the definition of MCP, $MCP(t)$ is between the minimum and maximum bid price of generators. The initial value for MCP is assumed to be the highest bid price of generators (i.e., $MCP_0(t)=Max\{B_i\}=111$). With the initial guess, the optimality subproblem (23) was solved, and then, the optimality cut was constructed and applied to the master problem. The master problem with the optimality cut was consequently solved to obtain new values for $y$ (i.e., $MCP(t)$). The second iteration starts with the new $\hat{Y}^1$ which is achieved from the master problem. Optimality subproblem based on new $\hat{Y}^1$ was infeasible; therefore, the feasibility subproblem was solved to create the first feasibility cut for the master problem. After 4 iterations, and after 2 feasibility cuts and 2 optimality cuts were applied to the master problem, solution convergence was achieved with the objective cost $1734015.125. The convergence of the GBD method has been presented in Fig. 3.

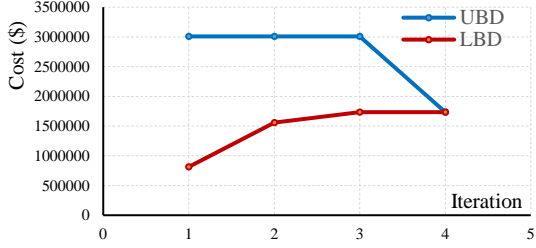

Fig. 3. Convergence of the proposed GBD-based algorithm

TABLE I
PAYMENT COST AND AVERAGE MCP FOR PCM AND OCM MECHANISMS WITH AND WITHOUT V2G

|  | OCM | OCM-V2G | PCM | PCM-V2G |
|---|---|---|---|---|
| Total Payment Cost ($) | 1823570 | 1888584.14 | 1816480 | 1734015.125 |
| Average MCP ($/MWh) | 63.125 | 66.167 | 62.792 | 60.417 |

The charging and discharging patterns of PEVs in OCM-V2G and PCM-V2G have been illustrated in Fig. 4(a). As shown in Fig. 4(a), PEVs are charged at the off-peak of demand profile, and they are discharged at the peak of load profile. As illustrated in Fig. 2(b) and Fig. 2(c), without V2G, PCM has the lower MCP and payment cost compared to those in OCM mechanism in the valley of load profile. However, when V2G is utilized in the system, PCM has a better performance compared to OCM not merely in the valley of load profile, it has also a better performance in other time periods, as shown in Fig. 4(b) and Fig. 4(c). It is because PCM mechanism uses the opportunity of charging and discharging PEVs to control and decrease MCPs and payment costs effectively. But, OCM mechanism uses charging and discharging PEVs to decrease the offer cost due to the objective function of OCM; therefore, because of the inconsistency between the clearing algorithm and the settlement rule, the MCP and payment cost in OCM-V2G are much higher than those in PCM-V2G mechanism. As presented in Table I, the average MCP of PCM-V2G (60.417 $/MWh), and also, the total payment cost of PCM-V2G ($1734015.125) are lower than the average MCP of OCM-V2G (66.167 $/MWh) and the total payment cost of OCM-V2G ($1888584.14), respectively. In fact, by using PCM instead of OCM in the system with V2G, $154569 is saved for this case study at the 24-hour period of time, and 8.18% saving is achieved. In summary, without V2G, PCM resulted in $7090 (0.38%) saving for payments compared to the OCM mechanism; however, this amount significantly increased to $154569 (8.18%) when V2G is utilized.

The curves of the net demand, MCP and the payment cost of OCM in the system with and without V2G are presented in Fig. 5. By using V2G in the system, in OCM mechanism, the offer cost is decreased from $1013048 to $993493.69; however, as shown in Fig. 5, MCP and payment cost have been increased due to the incompatibility between the clearing algorithm (i.e., OCM) and the settlement rule (i.e., pay-as-MCP). In fact, because of the objective function of OCM, OCM mechanism uses charging and discharging PEVs to decrease the offer cost without regarding to the increase of the payment costs and MCPs. Therefore, MCPs and payment costs in OCM in the system with V2G have been increased in comparison with those in the system without V2G.

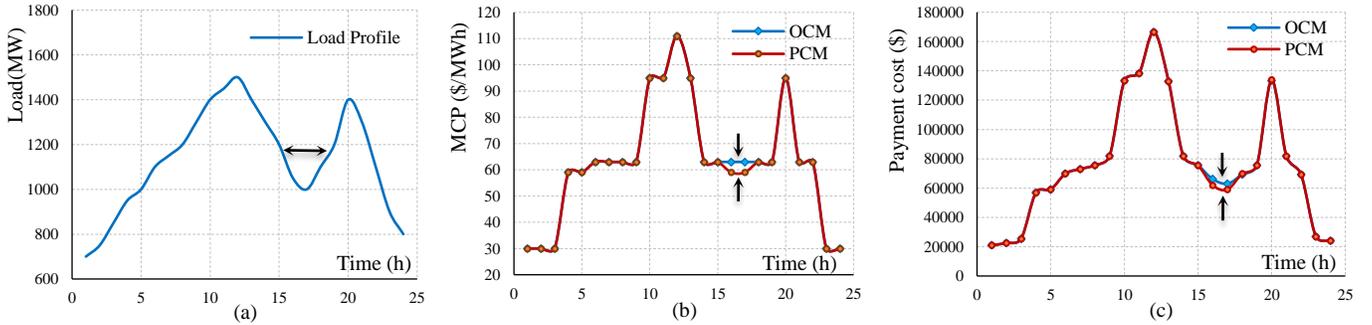

Fig. 2. Performance of PCM and OCM mechanisms in the system without V2G: (a) Load profile; (b) MCP; (c) Payment cost

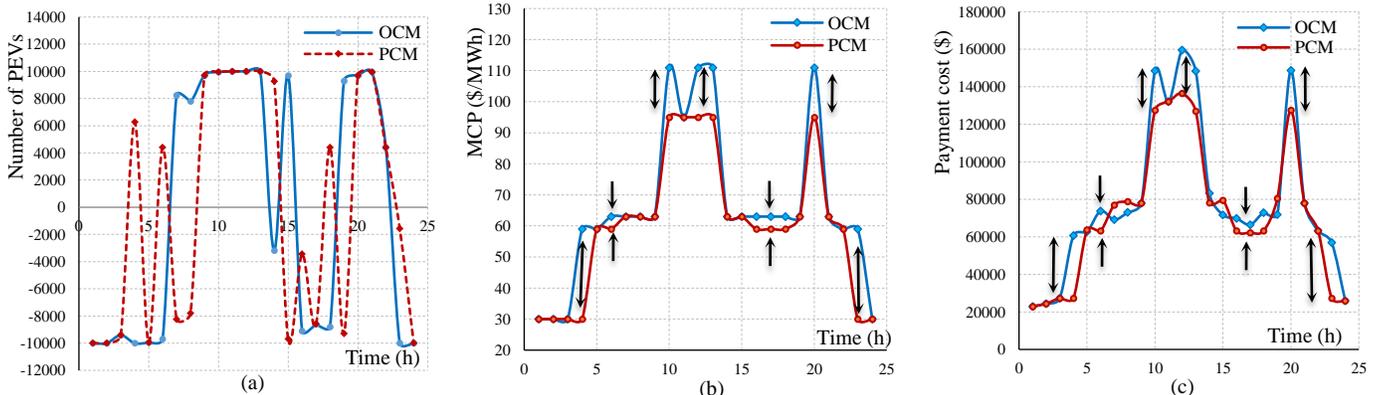

Fig. 4. Comparison of the performance of OCM and PCM mechanisms in the system with and without V2G:
(a) Charging and discharging pattern of PEVs; (b) MCP; (c) Payment cost

TABLE II
UNITS PAYMENTS IN PCM AND OCM MECHANISMS IN THE SYSTEM WITH AND WITHOUT V2G

| Unit | OCM | OCM-V2G | PCM | PCM-V2G |
|---|---|---|---|---|
| 1 | 687225 | 722540 | 585221 | 659750 |
| 2 | 674925 | 713877.5 | 669185 | 638955.5 |
| 3 | 113016 | 99895.31 | 122260 | 106031.1 |
| 4 | 99040 | 107192.11 | 146850 | 86050.47 |
| 5 | 198866 | 224533.71 | 220734 | 184624.3 |
| 6 | 29303 | 14373.75 | 39960 | 38680 |
| 7 | 20025 | 0 | 26105 | 19923.75 |
| 8 | 1170 | 1170 | 0 | 0 |
| 9 | 0 | 0 | 6165 | 0 |
| 10 | 0 | 5001.75 | 0 | 0 |

The curves of the net demand, MCP and the payment cost of PCM in the system with and without V2G are presented in Fig. 6. As shown in Fig. 6, by using V2G in the system, the performance of PCM has been improved, and MCPs and payment costs have been decreased due to the compatibility between the clearing algorithm (i.e., PCM) and the settlement rule (i.e., pay-as-MCP).

Payments to each unit for each OCM, OCM-V2G, PCM, and PCM-V2G problems are presented in Table II. As shown in Table II, in PCM mechanism, by using V2G in the system, unit 9 which had offered high price is not utilized, and payment cost decreases. In OCM mechanism, by using V2G, because of the goal of the OCM problem (i.e., minimizing offer cost), unit 10 with lower startup cost and higher bid price is selected instead of unit 7 with higher startup cost and lower bid price to decrease the offer cost, and it is the reason why despite decreasing the offer cost in OCM-V2G, the payment cost has increased.

To examine the performance of the proposed GBD-based method, the PCM-V2G problem was also solved by using two commercial solvers DICOPT [29], which is based on the extensions of the outer-approximation algorithm for the equality relaxation strategy, and SBB [30], which is based on a combination of the standard Branch and Bound method and some of the standard nonlinear programing software (NLP) solvers, as well as three free solvers SCIP [31], COUENNE [32], and BONMIN [33]. Table III demonstrates the effective performance of the proposed GBD-based solution method compared to MINLP solvers.

## V. CONCLUSIONS

This study has proposed a GBD-based solution method to solve PCM problems which is able to handle the nonlinearity posed by the V2G integration. Furthermore, in this study, the performances of PCM and OCM auction mechanisms has been

TABLE III
THE PERFORMANCE OF THE PROPOSED METHOD COMPARED TO OTHER METHODS FOR PCM-V2G PROBLEM

| Solution methods | Optimal Solution ($) | Convergence Time (Sec) | Relative Gap (%) |
|---|---|---|---|
| Proposed GBD-based method | 1734015.125 | 1.166 | 0 |
| DICOPT | 1734015.125 | 111.605 | 0 |
| SBB | 1972314.025 | 360 | 1.1581 |
| SCIP | 2020842.85 | 360 | 118.54 |
| COUENNE | 2132356.15 | 360 | 3.0119 |
| BONMIN | 2788054.25 | 360 | 1.8917 |

investigated in the system with and without V2G. The key finding in this paper are summarized as follows:

1) By using V2G in OCM mechanism, the offer cost were slightly decreased; however, MCP and payment cost were increased which is because of the incompatibility between the clearing algorithm (i.e., OCM) and the settlement rule (i.e., pay-as-MCP). But, by using V2G in the system, the performance of PCM were improved, and MCPs and payment costs were decreased due to the compatibility between the clearing algorithm (i.e., PCM) and the settlement rule (i.e., pay-as-MCP). Based on the numerical results, without using V2G, PCM resulted in $7090 (0.38%) saving for payments compared to the OCM mechanism; however, this amount significantly increased to $154569 (8.18%) when V2G is utilized. In conclusion, in the system in which V2G is utilized, using the PCM mechanism can considerably decrease the payment costs and achieve a huge saving compared to OCM mechanism.

2) It has been shown that the proposed GBD-based method is effective to solve the PCM problem which is a complex self-referring and NP-hard problem. Since the proposed GBD-based method split the PCM-V2G problem into simple MILP subproblems and LP master problem, it can be utilized in large scale systems with the increased number of decision variables and constraints. In addition, since the proposed method is an iterative method, linear sensitivity coefficients (LSC) can be also used to cope with the nonlinearity posed by nonlinear network constraints (either in transmission or distribution systems) as explained more in [34, 35]. In conclusion, the proposed GBD-based method presented in this paper can pave the way for future research works on PCM auction mechanism.

Future work may include nonlinear network constraints as well as the probabilistic characteristics of load variation, especially under the smart grid initiatives for demand response.

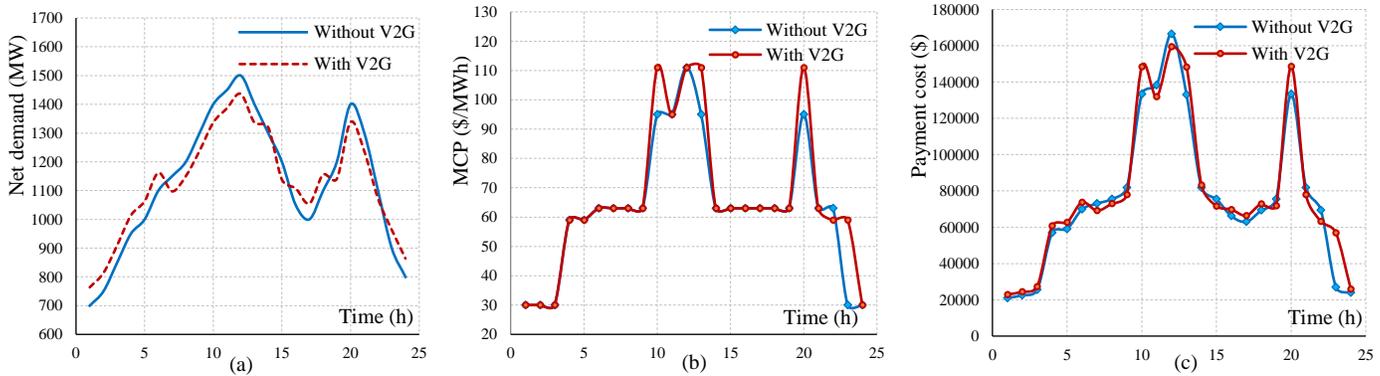

Fig. 5. Performance of OCM mechanism in the system with and without V2G: (a) Net demand; (b) MCP; (c) Payment cost

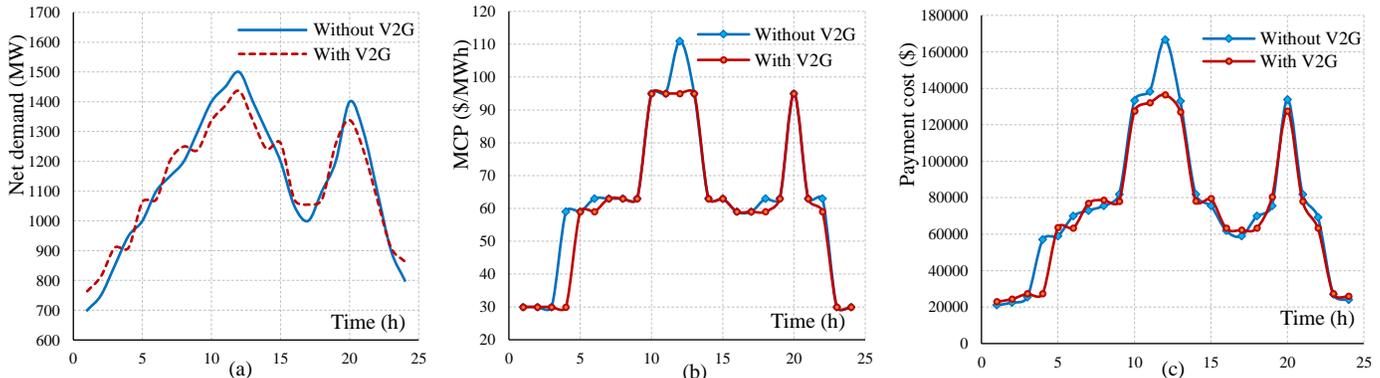

Fig. 6. Performance of PCM mechanism in the system with and without V2G: (a) Net demand; (b) MCP; (c) Payment cost